\documentclass[journal]{IEEEtran}
\ifCLASSINFOpdf
\else
\fi
\usepackage{cite}
\usepackage{amsmath}
\usepackage[pdftex]{graphicx} 
\usepackage{subfig}
\usepackage{array} 
\usepackage{fixltx2e}

\usepackage{booktabs} 
\usepackage{amssymb}
\usepackage{mathtools}
\usepackage{cancel} 
\usepackage{bm} 
\usepackage{esint} 
\usepackage{amsthm}
\usepackage{siunitx}
\usepackage{xfrac}
\usepackage{balance}

\usepackage{tikz}
\usepackage{tikz-3dplot}
\usetikzlibrary{backgrounds,3d,shadings,shapes.misc,decorations.pathmorphing,shapes,calc,fadings,shadows.blur}

\usepackage{amssymb,amsmath}

\usepackage{bm}
\newcommand{\bs}[1]{\bm{\mathrm{#1}}}

\newcommand{\vect}[1]{\vec{#1}}
\newcommand{\nhat}{\ensuremath{\hat{n}}}

\newcommand{\vecprime}[1]{\vect{#1}^{\,\prime}}

\renewcommand{\r}{\left(\vect{r}\right)}
\newcommand{\rp}{\left(\vecprime{r}\right)}

\newcommand{\rrpk}[1][]{\left(k_{#1}, \vect{r}, \vecprime{r}\right)}
\newcommand{\rrpki}[1][]{\left(k_i, \vect{r}, \vecprime{r}\right)}

\newcommand{\matr}[1]{\bs{#1}}

\newcommand{\abs}[1]{\left \lvert #1 \right\rvert }
 
\newcommand{\junk}[1] {}

\def\XXint#1#2#3{{\setbox0=\hbox{$#1{#2#3}{\int}$}
\vcenter{\hbox{$#2#3$}}\kern-.5\wd0}}

\newcommand*\widebar[1]{%
  \hbox{%
    \vbox{%
      \hrule height 0.5pt 
      \kern0.3ex
      \hbox{%
        \kern-0.05em
        \ensuremath{#1}%
        \kern-0.05em
      }%
    }%
  }%
} 

\newcommand{\mathsout}[1]
{\bgroup\mathchoice
  {\sbox0{$\displaystyle{#1}$}%
    \usebox0\hspace{-\wd0}%
    \rule[0.5\ht0-0.5\dp0-.5pt]{.7\wd0}{1pt}\hspace{.3\wd0}}%
  {\sbox0{$\textstyle{#1}$}%
    \usebox0\hspace{-\wd0}%
    \rule[0.5\ht0-0.5\dp0-.5pt]{.7\wd0}{1pt}\hspace{.3\wd0}}%
  {\sbox0{$\scriptstyle{#1}$}%
    \usebox0\hspace{-\wd0}%
    \rule[0.5\ht0-0.5\dp0-.5pt]{.7\wd0}{1pt}\hspace{.3\wd0}}%
  {\sbox0{$\scriptscriptstyle{#1}$}%
    \usebox0\hspace{-\wd0}%
    \rule[0.5\ht0-0.5\dp0-.5pt]{.7\wd0}{1pt}\hspace{.3\wd0}}%
\egroup}

\renewcommand{\epsilon}{\varepsilon}

\newcommand{\opL}{\ensuremath{\mathcal{L}}} 
\newcommand{\opK}{\ensuremath{\mathcal{K}}} 
\newcommand{\opM}{\ensuremath{\mathcal{M}}} 

\newcommand{\Lmat}[1][]{{\matr{L}_{#1}}}

\newcommand{\Kpvmat}[1][]{{\mathsout{\matr{K}}_{#1}}}

\newcommand{\Pmat}[1][]{{\matr{I}_{\parallel}^{#1}}}

\newcommand{\Mpvmat}[1][]{{\mathsout{\matr{M}}_{#1}}}

\renewcommand{\tt}{{(\mathrm{tt})}}

\newcommand{\Lttmat}[1][]{{\matr{L}_{#1}^{\mathrm{(ff)}}}}
\newcommand{\Lssmat}[1][]{{\matr{L}_{#1}^{\mathrm{(hh)}}}}

\newcommand{\Lnnmat}[1][]{{\matr{L}_{#1}^{\mathrm{(hh)}}}}
\newcommand{\Lntmat}[1][]{{\matr{L}_{#1}^{\mathrm{(hf)}}}}
\newcommand{\Ltnmat}[1][]{{\matr{L}_{#1}^{\mathrm{(fh)}}}}

\newcommand{\Mnspvmat}[1][]{{\mathsout{\matr{M}}_{#1}^{\mathrm{(hh)}}}}
\newcommand{\Msnpvmat}[1][]{{\mathsout{\matr{M}}_{#1}^{\mathrm{(hh)}}}}
\newcommand{\Kttpvmat}[1][]{{\mathsout{\matr{K}}_{#1}^{\mathrm{(fg)}}}}

\newcommand{\Kntpvmat}[1][]{{\mathsout{\matr{K}}_{#1}^{\mathrm{(hg)}}}}

\newcommand{\Irtmat}{{\matr{I}_{\times}^{\mathrm{(fg)}}}}

\newcommand{\figref}[1]{Fig.~\ref{#1}}
\newcommand{\secref}[1]{Section~\ref{#1}}

\usepackage{tcolorbox}

\usepackage{textpos}

\newcommand{\mySubtitle}[1]%
{%
	\begin{textblock}{14.0}(0.7, 2.9)
		\textbf{#1}%
	\end{textblock}%
}%

\newcommand{\green}[1]{\textcolor{green!0!black}{#1}}

\newcommand{\brown}[1]{\textcolor{black}{#1}}
\newcommand{\magenta}[1]{\textcolor{black}{#1}}

\newcommand{\Ecolor}[1]{{#1}}
\newcommand{\Hcolor}[1]{{#1}}
\newcommand{\Jcolor}[1]{{#1}}

\newcommand{\Rhocolor}[1]{\brown{#1}}

\newcommand{\Er}[1][]{\Ecolor{\vect{E}_{#1}\r}}

\newcommand{\Etr}[1][]{\Ecolor{\nhat_{#1} \times \vect{E}_{#1}\r}}

\newcommand{\Hr}[1][]{\Hcolor{\vect{H}_{#1}\r}}

\newcommand{\Htr}[1][]{\Hcolor{\nhat_{#1} \times \vect{H}_{#1}\r}}

\newcommand{\Jr}[1][]{\ensuremath{\Jcolor{\vect{J}_{#1}\r}}}

\newcommand{\Grrpk}[1][]{\ensuremath{\green{G_{#1}\rrpk}}}

\newcommand{\Grrpki}[1][]{\ensuremath{\green{G\rrpki}}}

\newcommand{\gradrGrrpk}[1][]{\ensuremath{\green{\nabla G_{#1}\rrpk}}}

\newcommand{\Ar}[1][]{\magenta{\vect{A}_{#1}\r}}

\newcommand{\Atr}[1][]{\magenta{\nhat \times \vect{A}_{#1}\r}}
\newcommand{\Atrp}[1][]{\magenta{\nhat' \times \vect{A}_{#1}\rp}}

\newcommand{\Anr}[1][]{\magenta{\nhat \cdot \vect{A}_{#1}\r}}
\newcommand{\Anrp}[1][]{\magenta{\nhat' \cdot \vect{A}_{#1}\rp}}

\newcommand{\divrAr}[1][]{\magenta{\nabla \cdot \vect{A}_{#1}\r}}

\newcommand{\curlrAr}[1][]{\magenta{\nabla \times \vect{A}_{#1}\r}}
\newcommand{\curlrpArp}[1][]{\magenta{\nabla' \times \vect{A}_{#1}\rp}}

\newcommand{\Phimat}[1][]{\magenta{\matr{\Phi}_{#1}}}

\newcommand{\phir}[1][]{\magenta{\phi_{#1}\r}}
\newcommand{\phirp}[1][]{\magenta{\phi_{#1}\rp}}

\newcommand{\gradrphir}[1][]{\magenta{\nabla\phi_{#1}\r}}
\newcommand{\gradrpphirp}[1][]{\magenta{\nabla\phi_{#1}\rp}}

\newcommand{\uamat}[1][]{\ensuremath{\Jcolor{\matr{a}_{c,#1}}}}
\newcommand{\ubmat}[1][]{\ensuremath{\Ecolor{\matr{a}_{t,#1}}}}

\newcommand{\udmat}[1][]{\ensuremath{\Rhocolor{\matr{a}_{n,#1}}}}

\newcommand{\ndgPhimat}[1][]{\ensuremath{\Rhocolor{\matr{\Psi}_{#1}}}}

\makeatletter
\newlength\numerator@height
\newlength\frac@height
\newsavebox\numerator@box
\newsavebox\frac@box

\newcommand\dfracparens[3]{%
	\sbox{\numerator@box}{\ensuremath{#1}}%
	\sbox{\frac@box}{\ensuremath{\dfrac{#1}{#2}}}%
	\settoheight{\frac@height}{\usebox{\frac@box}}%
	\settoheight{\numerator@height}{\usebox{\numerator@box}}%
	\addtolength{\frac@height}{-\numerator@height}%
	\usebox{\frac@box}%
	\raisebox{\frac@height}{%
		\( \left( {#3} \right)
		\)}%
}
\makeatother

\begin{document}
%
%
\title{Electromagnetic Modeling of Lossy~Materials with a Potential-Based Boundary~Element~Method}

%
%
%

\author{Shashwat~Sharma,~\IEEEmembership{Graduate Student Member,~IEEE,}
        and~Piero~Triverio,~\IEEEmembership{Senior Member,~IEEE}
\thanks{S. Sharma is with the Edward S. Rogers Sr. Department of Electrical \& Computer Engineering, University of Toronto, Toronto, ON, M5S 3G4 Canada, e-mail: shash.sharma@mail.utoronto.ca.
	P. Triverio is with the Edward S. Rogers Sr. Department of Electrical \& Computer Engineering and with the Institute of Biomedical Engineering, University of Toronto, Toronto, ON, M5S 3G4 Canada, email: piero.triverio@utoronto.ca.}
\thanks{This work was supported by Advanced Micro Devices, by the Natural Sciences and Engineering Research Council of Canada (Collaborative Research and Development Grants program), and by CMC Microsystems.}
\thanks{Manuscript received $\ldots$; revised $\ldots$.}}

%
%

\markboth{IEEE Antennas and Wireless Propagation Letters}%
{Sharma \MakeLowercase{\textit{et al.}}: Potential-Based EM Modeling of Lossy Materials}
%



\maketitle

\begin{abstract}
The boundary element method (BEM) enables solving three-dimensional electromagnetic problems using a two-dimensional surface mesh, making it appealing for applications ranging from electrical interconnect analysis to the design of metasurfaces. The BEM typically involves the electric and magnetic fields as unknown quantities. Formulations based on electromagnetic potentials rather than fields have garnered interest recently, for two main reasons: (a) they are inherently stable at low frequencies, unlike many field-based approaches, and (b) potentials provide a more direct interface to quantum physical phenomena. Existing potential-based formulations for electromagnetic scattering have been proposed primarily for perfect conductors. We develop a potential-based BEM formulation which can capture both dielectric and conductive losses, and accurately models the skin effect over broad ranges of frequency. The accuracy of the proposed formulation is validated through canonical and realistic numerical examples. 
\end{abstract}


%
\IEEEpeerreviewmaketitle


\section{Introduction}

\IEEEPARstart{E}{lectromagnetic} simulation tools based on the boundary element method (BEM)~\cite{ChewWAF} have gained traction in a variety of applications, ranging from antenna modeling~\cite{PMCHWT01} to high-speed interconnect analysis~\cite{gibc,gibcHmatDanJiao,DSA08,TAPAIMin,AWPLSLIM}.
The BEM is based on a surface integral representation of Maxwell's equations, which allows three-dimensional problems to be solved in terms of quantities defined on a two-dimensional surface mesh.
Many conventional full-wave BEM formulations, which take electric and magnetic fields as unknown quantities, suffer from numerical instability at very low frequencies~\cite{lfbreakdown}.
The need to model multiscale structures at both high and low frequencies arises in applications such as the analysis of integrated circuit components.
This motivates the development of broadband BEM formulations where electromagnetic potentials, rather than fields, are taken as the unknowns~\cite{PIE01,PIE04}.
These formulations do not rely on the coupling between electric and magnetic fields, and provide a natural interface to quantum phenomena, which is an important consideration in emerging applications such as quantum computing~\cite{Chew_QEM1_JMMCT,Chew_QEM2_JMMCT}.
Potential-based methods may also be well-suited for coupled electromagnetic-circuit simulations.

Potential-based integral equation (PIE) methods for scattering analysis have been developed primarily for perfect conductors~\cite{PIE01,PIE04,PIE03,PIE08,PIE_TD_PEC}.
For low frequencies and sub-wavelength structures, where PIE formulations are most sorely needed, the skin depth in a conductor may be large compared to its physical dimensions, and modeling the conductor as perfect may be inaccurate.
Although magnetoquasistatic PIE formulations have been used for eddy current modeling in lossy conductors~\cite{eddy01,eddy02,eddy03,eddy04,eddy05}, these techniques do not apply in the presence of dielectric inclusions, or at high frequencies.
PIE-based modeling of lossless dielectric objects was considered in~\cite{PIE02}, but mainly from a theoretical perspective, and in~\cite{PIE_TD_DIEL}, but in the time domain.
Also, the method in~\cite{PIE02} involves a linear combination of integral equations written for adjacent materials.
As in the case of analogous field-based formulations~\cite{PMCHWT02,PMCHWT03,PMCHWT04}, the formulation in~\cite{PIE02} may be inaccurate for large contrasts in material parameters of adjacent media.
To the best of our knowledge, a full-wave PIE formulation for lossy penetrable materials has not been demonstrated in the frequency domain.
 
In this article, we devise a novel full-wave PIE formulation for electromagnetic scattering from lossy dielectrics and conductors, applicable at both low and high frequencies.
The proposed formulation couples the scalar and vector potential integral equations~\cite{PIE01} in the regions internal and external to each object.
An appropriate discretization is discussed, and the accuracy of the formulation is demonstrated numerically over wide ranges of conductivity and frequency.

\section{Proposed Formulation}

We consider time-harmonic scattering from an object occupying volume~$\mathcal{V}$, with surface~$\mathcal{S}$ and outward unit normal vector~$\nhat$.
Symbols~$\mathcal{S}^-$ and~$\mathcal{S}^+$ denote the internal and external sides of~$\mathcal{S}$, respectively.
The object has permittivity~$\epsilon$, permeability~$\mu$, and conductivity~${\sigma > 0}$.
The permittivity may be complex,~${\epsilon = \epsilon' - j\epsilon''}$, where the imaginary part represents dielectric losses~\cite{pozar}.
The object resides in free space,~$\mathcal{V}_0$, with permittivity~$\epsilon_0$ and permeability~$\mu_0$.

\subsection{Internal Region}\label{sec:int}

For~${\vect{r}\in\mathcal{V}}$, the magnetic vector and electric scalar potentials,~$\Ar$ and~$\phir$, respectively, can be defined via~\cite{jacksonEM}
\begin{align}
	\mu\Hr &= \nabla\times\Ar, \label{Adef}\\
	\Er &= -j\omega\Ar - \gradrphir, \label{phidef}
\end{align}
where~$\omega$ is the angular frequency,~$\Hr$ is the magnetic field, and~$\Er$ is the electric field.
Using~\eqref{Adef} and~\eqref{phidef} in Maxwell's equations, it can be shown that~$\Ar$ satisfies the homogeneous Helmholtz equation~\cite{jacksonEM}
\begin{align}
	\nabla^2\Ar + k^2\Ar = 0, \label{Ahelm}
\end{align}
when the Lorenz gauge
\begin{align}
	\divrAr = -\mu\left(j\omega\epsilon' + \omega\epsilon'' + \sigma\right)\,\phir \label{lorenz}
\end{align}
is adopted~\cite{rothwellEM,PIE01}.
In~\eqref{Ahelm},~${k^2=-j\omega\mu\left(j\omega\epsilon' + \omega\epsilon'' +\sigma\right)}$ is the wave number associated with the object's material.
Likewise,~$\phir$ satisfies the Helmholtz equation~\cite{jacksonEM}
\begin{align}
	\left(\gamma\mu^{-1} + \sigma\right)\left[\nabla^2\phir + k^2\phir\right] = 0,\label{phihelm}
\end{align}
where~${\gamma \triangleq \left(j\omega\epsilon' + \omega\epsilon'' + \sigma\right)\mu}$.
The use of Green's identities~\cite{HansonYakovlev} with~\eqref{Ahelm} and~\eqref{phihelm} leads to the integral equations~\cite{PIE01}
\begin{multline}
    \opL\bigl[\nhat'\times\curlrpArp\bigr] + \opK\bigl[\Atrp\bigr]\\
    + \Ar
    + \gamma\,\opL\bigl[\phirp\nhat'\bigr]\\
    -\nabla\opL\bigl[\Anrp\bigr] = 0,\label{vpieint}
\end{multline}
\begin{align}
  \gamma\mu^{-1}\left(\opL\bigl[\nhat'\cdot\gradrpphirp\bigr]
  + \opM\bigl[\phirp\bigr] - \phir\right)
  = 0,\label{spieint}
\end{align}
where primed and unprimed coordinates denote source and observation points, respectively,~${\vect{r}\in\mathcal{V}}$, and~${\vecprime{r}\in\mathcal{S}^-}$.
The integral operators in~\eqref{vpieint} and~\eqref{spieint} are defined as~\cite{ChewWAF}
\begin{align}
	\opL\bigl[\vect{a}\rp\bigr] &= \int_{\mathcal{S}}d\mathcal{S}'\,\Grrpk\,\vect{a}\rp,\label{opLdef}\\
	\opK\bigl[\vect{a}\rp\bigr] &= \int_\mathcal{S}d\mathcal{S}'\,\nabla\Grrpk\times\vect{a}\rp,\label{opKdef}\\
	\opM\bigl[a\rp\bigr] &= \int_\mathcal{S}d\mathcal{S}'\,\nhat'\cdot\gradrGrrpk\,a\rp,\label{opMdef}
\end{align}
where~$\Grrpk$ is the Green's function associated with the object's material,
\begin{align}
	\Grrpk = \frac{e^{-jk\abs{\vect{r}-\vecprime{r}}}}{4\pi\abs{\vect{r}-\vecprime{r}}}.\label{grrpk}
\end{align}

Next, the object's surface is discretized with a triangular mesh.
Quantity~${\nhat\times\curlrAr}$ is expanded with Rao-Wilton-Glisson (RWG) functions~\cite{RWG},~$\vect{f}_n\r$ normalized by edge length, while~$\Atr$ is expanded with Buffa-Christiansen functions~\cite{BCorig},~$\vect{g}_n\r$, which are defined on a barycentric refinement of the mesh.
This choice of functions stems from~\eqref{Adef} and~\eqref{phidef}, which indicate that~$\nhat\times\curlrAr$ and~$\Atr$ are related to~$\Htr$ and~$\Etr$ on~$\mathcal{S}$, respectively.
Since~$\Htr$ and~$\Etr$ can be interpreted as electric and magnetic surface current densities, respectively,~$\nhat\times\curlrAr$ and~$\Atr$ must both be expanded with divergence-conforming basis functions~\cite{ChewIEM}, while respecting their mutual orthogonality~\cite{calderon01}.
These requirements are satisfied by the proposed expansion scheme.

Unknowns~$\phir$ and~$\Anr$ are expanded with unit-amplitude pulse functions,~$h_n\r$, while~$\nhat\cdot\gradrphir$ is expanded with area-normalized pulse functions,~$\widetilde{h}_n\r$.
The choice of normalizing by area was based on an empirical study of the condition number of the final system matrix.
A more sophisticated choice of basis function for~$\phir$, such as the one suggested in~\cite{PIE_TD_DIEL}, may lead to improved accuracy.
However, as shown in~\secref{sec:results}, the pulse functions are sufficiently accurate in a variety of cases.

Taking the cross product of~$\nhat$ with~\eqref{vpieint} and letting~$\vect{r}\to\mathcal{S}^-$ gives its rotated tangential part, which is tested with~$\nhat\times\text{RWG}$ functions to get the matrix relation
\begin{multline}
	\Lttmat\uamat + \left(\Kttpvmat + \frac{1}{2}\,\Irtmat\right)\ubmat\\ + \gamma\,\Ltnmat\Phimat + \matr{D}^T\Lssmat\udmat = \matr{0}.\label{tvpieintd}
\end{multline}
Taking the dot product of~$\nhat$ with~\eqref{vpieint} and letting~$\vect{r}\to\mathcal{S}^-$ gives its normal component, which is tested with $\widetilde{h}_n\r$ to obtain
\begin{multline}
	\Lntmat\uamat + \Kntpvmat\ubmat + \gamma\,\Lnnmat\Phimat\\ - \left(-\left(\Mnspvmat\right)^T - \frac{1}{2}\Pmat[(\mathrm{hh})]\right)\udmat = \matr{0},\label{nvpieintd}
\end{multline}
where the superscript~``$T$'' denotes the matrix transpose.
Finally, taking~$\vect{r}\to\mathcal{S}^-$ in~\eqref{spieint} and testing it with $\widetilde{h}_n\r$ yields
\begin{align}
	\gamma\mu^{-1}\left[\Lssmat\ndgPhimat + \left(\Mnspvmat - \frac{1}{2}\Pmat[(\mathrm{hh})]\right)\Phimat\right] = \matr{0}.\label{spieintd}
\end{align}
In~\eqref{tvpieintd},~\eqref{nvpieintd}, and~\eqref{spieintd},~$\Lmat$,~$\Kpvmat$ and~$\Mpvmat$ are the discretized~$\opL$,~$\opK$ and~$\opM$ operators, respectively, where a dash through a matrix indicates that the associated integral is computed in the principal value sense~\cite{ChewWAF}. Term~$\matr{D}$ is a sparse incidence matrix linking mesh edges and triangles~\cite{aefie2}. Identity operator~$\Irtmat$ is obtained by testing~$\vect{g}_n\r$ with~$\nhat\times\vect{f}_n\r$, while~$\Pmat[(\mathrm{hh})]$ involves testing $h_n\r$ with $\widetilde{h}_n\r$.
The superscript labels~$(ij)$ on each discrete operator represent the testing and basis functions involved, respectively.
Column vectors~$\uamat$,~$\ubmat$,~$\udmat$,~$\Phimat$ and~$\ndgPhimat$ contain the unknown coefficients associated with~${\nhat\times\curlrAr}$,~${\Atr}$,~${\Anr}$,~${\phir}$ and~${\nhat\cdot\gradrphir}$, respectively.

\subsection{External Region}\label{sec:ext}

Next, PIEs are devised to capture the physics in the region external to~$\mathcal{S}$, following the procedure in \secref{sec:int} for~$\vect{r},\vecprime{r}\in\mathcal{S}^+$.
However, instead of using the normal component of~\eqref{vpieint} for the external region, we take its divergence~\cite{PIE03} and test the resulting equation with~$\widetilde{h}_n\r$.
We found that this choice leads to better conditioning of the final system matrix, and better accuracy.
The resulting discrete equations are
\begin{multline}
	\Lttmat[0]\uamat[0] + \left(\Kttpvmat[0] - \frac{1}{2}\,\Irtmat\right)\ubmat[0]\\ + \gamma_0\,\Ltnmat[0]\Phimat[0] + \matr{D}^T\Lssmat[0]\udmat[0] = -\ubmat[\mathrm{inc}],\label{tvpieextd}
\end{multline}
\begin{multline}
	\Lssmat[0]\matr{D}\uamat[0]
	+ \gamma_0\left(\Msnpvmat[0] + \frac{1}{2}\,\Pmat[(\mathrm{hh})]\right)\Phimat[0]\\
	+ k_0^2\,\Lssmat[0]\udmat[0]
	= \gamma_0\,\Phimat[\mathrm{inc}],\label{nvpieextd}
\end{multline}
\begin{align}
	\Lssmat[0]\ndgPhimat[0] + \left(\Mnspvmat[0] + \frac{1}{2}\Pmat[(\mathrm{hh})]\right)\Phimat[0] = -\Phimat[\mathrm{inc}],\label{spieextd}
\end{align}
where the subscript~``$0$'' on the matrix operators denotes that the Green's function associated with~$\mathcal{V}_0$ is used, and~${\gamma_0 \triangleq j\omega\epsilon_0\mu_0}$.
Subscript~``$0$'' on the column vectors of unknowns indicates that the quantities are defined on~$\mathcal{S}^+$.
Subscript~``$\mathrm{inc}$'' indicates incident potentials~\cite{PIE01,PIE04}.

\subsection{Boundary Conditions}\label{sec:bc}

For~${\nhat\times\curlrAr}$,~${\Atr}$ and~${\phir}$, we use boundary conditions on~$\mathcal{S}$ identical to those in~\cite{PIE01,PIE04},
\begin{align}
	\mu_0^{-1}\,\nhat\times\curlrAr[0] &= \mu^{-1}\,\nhat\times\curlrAr,\label{bcAc}\\
	\Atr[0] &= \Atr,\label{bcAt}\\
	\phir[0] &= \phir,\label{bcPhi}
\end{align}
where~$\Ar$ and~$\Ar[0]$ are the magnetic vector potentials on~$\mathcal{S}^-$ and~$\mathcal{S}^+$, respectively.
For conductive objects, a new boundary condition is required for~${\Anr}$ and~${\nhat\cdot\gradrphir}$,
\begin{multline}
	-j\omega\left[\gamma_0\mu_0^{-1}\Anr[0] - \gamma\mu^{-1}\Anr\right]\\ - \left[\gamma_0\mu_0^{-1}\nhat\cdot\gradrphir[0] - \gamma\mu^{-1}\nhat\cdot\gradrphir\right] = 0,\label{bcAn}
\end{multline}
which is used to eliminate~$\nhat\cdot\gradrphir$ in~\eqref{spieint}.
Equation~\eqref{bcAn} is derived using~\eqref{phidef},~\eqref{lorenz} and the standard boundary conditions~\cite{rothwellEM} for~$\nhat\cdot\Er$,~$\nhat\cdot\Er[0]$ and ~$\nhat\cdot\Jr$, where~$\Er$ and~$\Er[0]$ are the electric fields in~$\mathcal{V}$ and~$\mathcal{V}_0$, respectively, and~$\Jr$ is the conduction volume current density in~$\mathcal{V}$.

\subsection{Final System Matrix}\label{sec:sys}

Equations~\eqref{tvpieextd},~\eqref{tvpieintd},~\eqref{nvpieextd},~\eqref{nvpieintd},~\eqref{spieintd} and~\eqref{spieextd} are concatenated in that order, and the boundary condition~\eqref{bcAn} is applied to get the final system of equations~\eqref{sys} at the top of the following page,
\begin{figure*}[t]
\small
\setcounter{equation}{21}
\begin{align}
	\setlength{\arraycolsep}{3pt}
    \begin{bmatrix}
		\frac{1}{\xi}\Lttmat[0] & \Kttpvmat[0] - \frac{1}{2}\Irtmat & \frac{jk_0}{\xi}\,\Ltnmat[0] & \matr{0} & \frac{1}{\xi}\matr{D}^T\Lssmat[0] & \matr{0} \\
		\frac{\mu}{\xi\mu_0}\Lttmat & \Kttpvmat + \frac{1}{2}\Irtmat & \frac{c_0\gamma}{\xi}\,\Ltnmat & \frac{1}{\xi}\matr{D}^T\Lssmat & \matr{0} & \matr{0} \\
		\Lssmat[0]\matr{D} & \matr{0} & jk_0\left(\Msnpvmat[0] + \frac{1}{2}\Pmat[(\mathrm{hh})]\right) & \matr{0} & k_0^2\,\Lssmat[0] & \matr{0} \\
		\frac{\mu}{\mu_0}\Lntmat & \xi\,\Kntpvmat & c_0\gamma\,\Lnnmat & \bigl(\Mnspvmat\bigr)^T - \frac{1}{2}\Pmat[(\mathrm{hh})] & \matr{0} & \matr{0} \\
		\matr{0} & \matr{0} & \Mnspvmat - \frac{1}{2}\Pmat[(\mathrm{hh})] & \frac{k^2}{\gamma c_0}\,\Lssmat & -\frac{k_0^2}{\gamma c_0}\frac{\mu}{\mu_0}\,\Lssmat & \frac{\gamma_0}{\gamma c_0}\frac{\mu}{\mu_0}\,\Lssmat \\
		\matr{0} & \matr{0} & \Mnspvmat[0] - \frac{1}{2}\Pmat[(\mathrm{hh})] & \matr{0} & \matr{0} & \frac{1}{c_0}\Lssmat[0]
	\end{bmatrix}
	{\renewcommand*{\arraystretch}{1.4}
	\begin{bmatrix}
		\uamat[0] \\ \ubmat/\xi \\ \Phimat/c_0 \\ \udmat \\ \udmat[0] \\ \ndgPhimat[0]
	\end{bmatrix}}
	=
	{\renewcommand*{\arraystretch}{1.4}
	\begin{bmatrix}
		-\ubmat[\mathrm{inc}]/\xi \\ \matr{0} \\ \gamma_0\,\Phimat[\mathrm{inc}] \\ \matr{0} \\ \matr{0} \\ -\Phimat[\mathrm{inc}]/c_0
	\end{bmatrix}}.%
	\label{sys}
\end{align}
\hrulefill
\end{figure*}
where~$\xi$ is the average mesh edge length.
In~\eqref{sys}, the equations and unknowns have been strategically scaled to ensure stable direct factorization and good accuracy for wide ranges of frequency and conductivity.
It may be necessary to enforce charge conservation at extremely low frequencies.
Although this point is not addressed in the present formulation, the results in \secref{sec:results} demonstrate the extremely wide ranges of frequency over which the proposed method still remains accurate.
Also, the system matrix in~\eqref{sys} contains two vector quantities ($\uamat[0]$,~$\ubmat$) and four scalar quantities ($\Phimat$,~$\udmat$,~$\udmat[0]$~$\ndgPhimat[0]$) as unknowns, while most field-based formulations for lossy conductors contain only two vector and up to one scalar unknown quantity~\cite{gibc,AWPLSLIM,eaefie02}.
However, the additional unknowns in the proposed method may be well worth the broadband performance, particularly when coupled with acceleration algorithms~\cite{MLFMA,AIMbles,pfftmain,TAPAIMx}.
It may be possible to obtain a smaller system matrix, but the scope of this work is to establish that potential-based formulations can be used to model lossy conductors accurately and over very wide ranges of frequency and conductivity.

\section{Results}\label{sec:results}

The accuracy of the proposed formulation is validated through comparisons with analytical solutions and an existing field-based BEM formulation: the enhanced augmented electric field integral equation~(eAEFIE) for penetrable objects~\cite{eaefie01,eaefie02}.
The numerical integration routines in~\cite{gibc} were used for the operators associated with the internal region, to maintain accuracy for highly conductive media.
Electric and magnetic fields tangential to~$\mathcal{S}$ were obtained as a post-processing step via~$\eqref{Adef}$ and~$\eqref{phidef}$.

\subsection{Sphere}\label{sec:results:sph}

\begin{figure}[t]
	\centering
	\includegraphics[width=\linewidth]{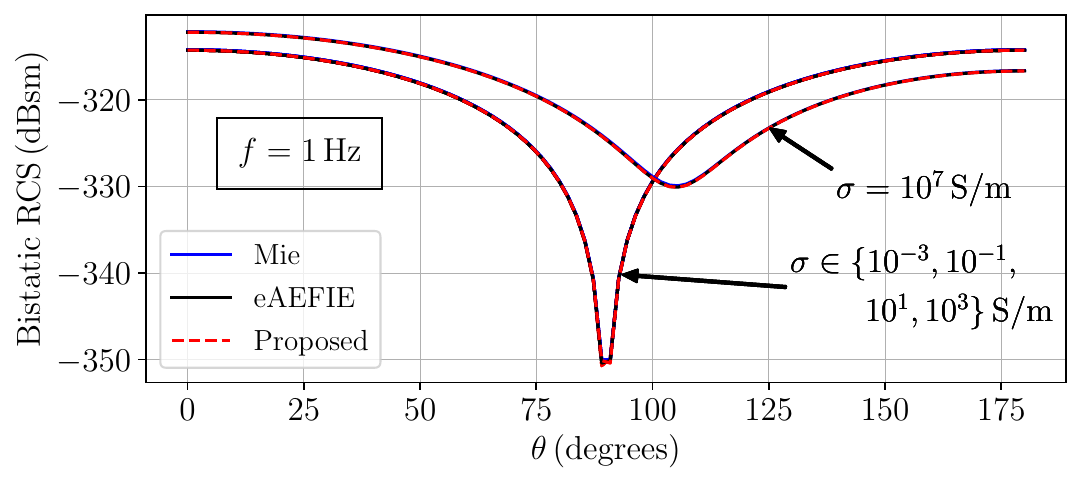}\\
	\vspace{-7mm}
	\includegraphics[width=\linewidth]{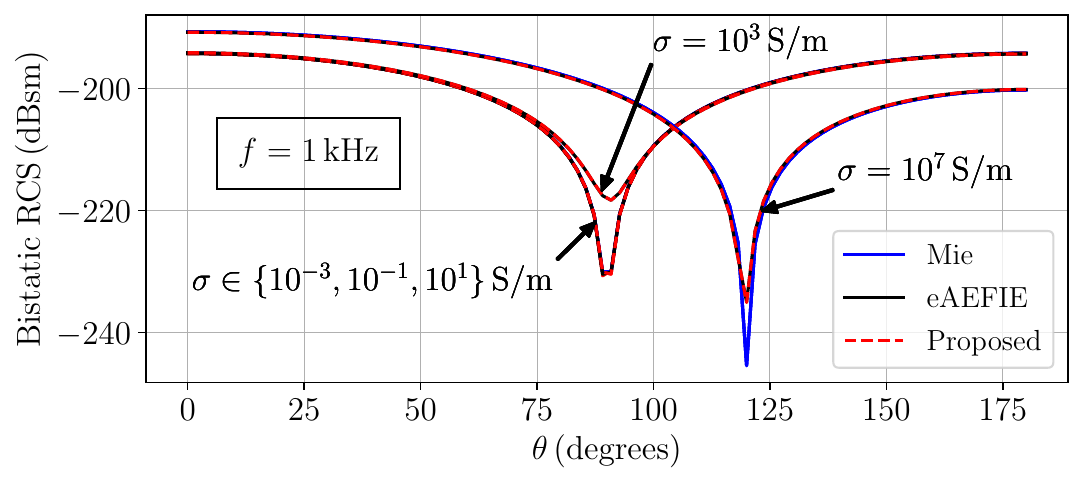}\\
	\vspace{-7mm}
	\includegraphics[width=\linewidth]{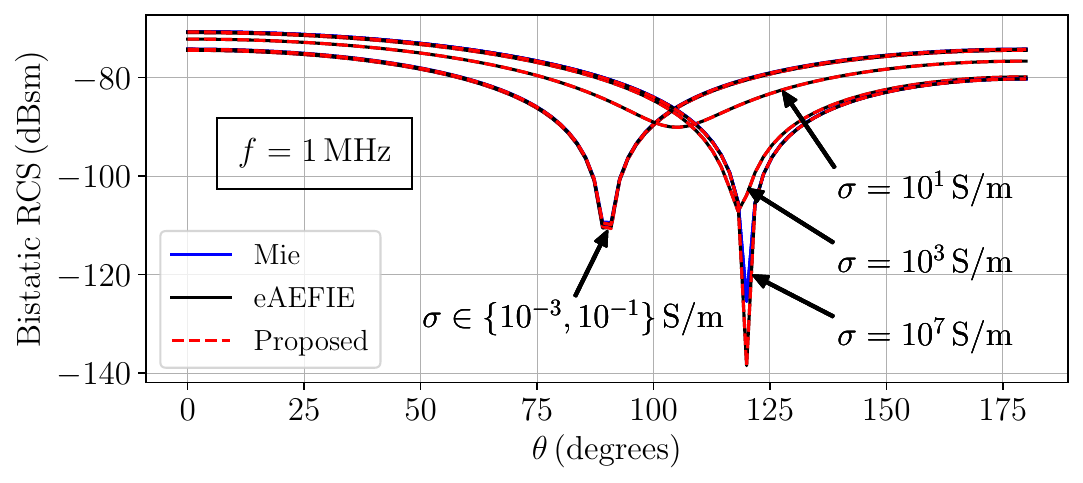}\\
	\vspace{-7mm}\hspace*{0.618mm}
	\includegraphics[width=0.983\linewidth]{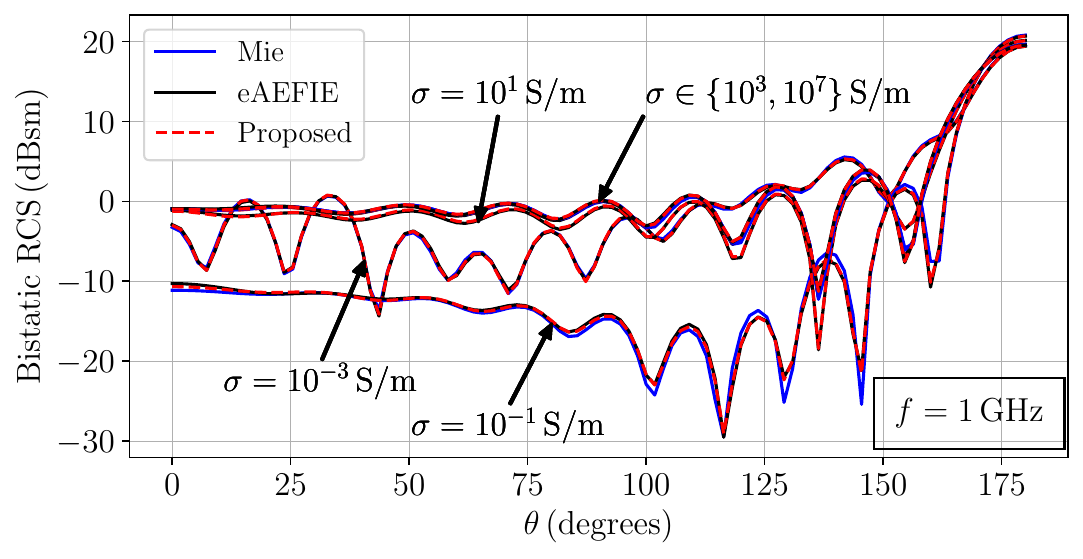}
	\caption{Accuracy validation for the sphere in \secref{sec:results:sph}. The same horizontal axis applies to all panels.}\label{fig:sph}
\end{figure}

First, we consider a sphere with diameter~$1\,$m and relative permittivity~$2$, excited by a plane wave.
The sphere is meshed with~$2,114$ triangles, and the bistatic radar cross section (RCS) is compared against the analytical Mie series (\figref{fig:sph}). The RCS is reported for the plane along which the incident electric field is polarized ($E$-plane).
We consider conductivities spanning~$10$ decades from~$10^{-3}\,$S/m, where the sphere behaves like a dielectric, to~$10^7\,$S/m, corresponding to a good conductor.
Nine decades of frequencies from~$1\,$Hz to~$1\,$GHz are simulated to encompass skin depths from~\SI{1}{\micro\metre} to several times the sphere's diameter.
For the~$1\,$GHz cases, a finer mesh with~$3,786$ triangles was used.
\figref{fig:sph} demonstrates the excellent accuracy of the proposed formulation.
The results deviate slightly from the Mie series in the bottom panel, but identical deviations are observed for the eAEFIE formulation.
Therefore, these errors can be attributed to the discretization and numerical integration, which are common to both formulations.

\subsection{Cube}\label{sec:results:cube}

\begin{figure}[t]
	\centering
	\includegraphics[width=\linewidth]{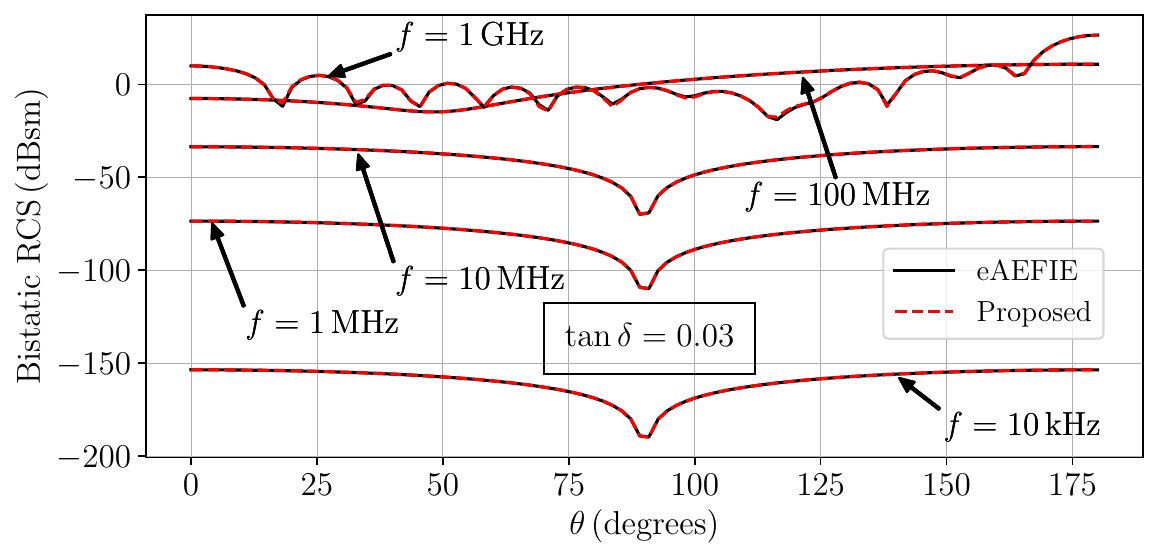}\\
	\vspace{-6.55mm}\hspace*{0.025mm}
	\includegraphics[width=.986\linewidth]{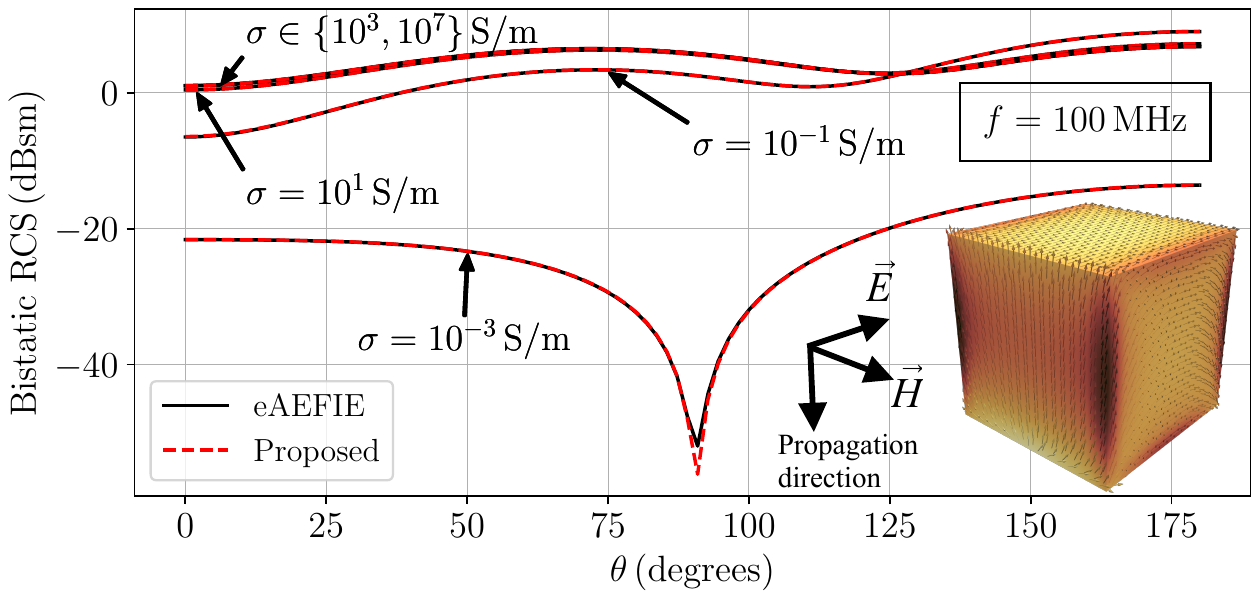}
	\caption{Accuracy validation for the cube in \secref{sec:results:cube}.}\label{fig:cube}
\end{figure}

Next, we consider an FR-4 cube meshed with~$3,274$ triangles, with side length~$1\,$m,~${\epsilon'=4.4\,\epsilon_0}$, and loss tangent~${\tan\delta=\sfrac{\epsilon''}{\epsilon'}=0.03}$~\cite{FR4}.
A plane wave impinges on the cube with the electric field polarized as shown in the inset in the bottom panel of \figref{fig:cube}.
The $E$-plane bistatic RCS for the proposed formulation is compared with the results obtained via the eAEFIE, for frequencies between~$10\,$kHz and~$1\,$GHz.
The top panel of \figref{fig:cube} demonstrates the excellent accuracy of the proposed formulation compared to the eAEFIE.
We also considered the case when~$\epsilon'=\epsilon_0$,~$\epsilon''=0$, and~$\sigma\in\left[10^{-3},10^7\right]\,$S/m at~$100\,$MHz, corresponding to skin depths between~\SI{1.6}{\metre} and~\SI{16}{\micro\metre}. Excellent agreement with the eAEFIE was achieved, as shown in the bottom panel of \figref{fig:cube}.
The inset shows the surface current density for~${\sigma=10}\,$S/m.

\subsection{Split Ring Resonator Array}\label{sec:results:srr}

\begin{figure}[t]
	\centering
	\includegraphics[width=\linewidth]{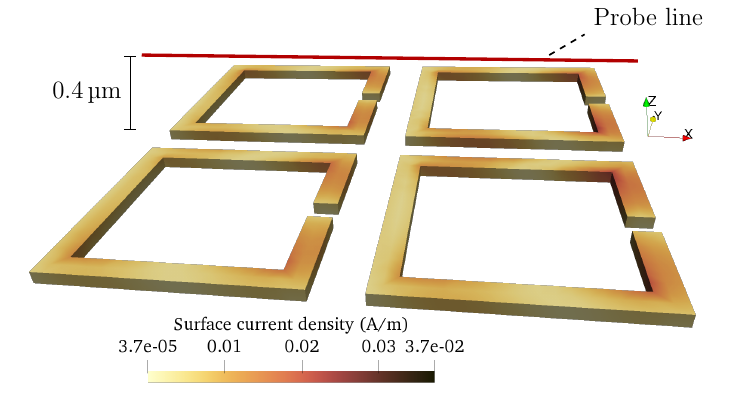}
	\caption{Geometry and electric surface current density magnitude for the SRR array in \secref{sec:results:srr}.}\label{fig:srrgeom}\vspace{1mm}
	\includegraphics[width=\linewidth]{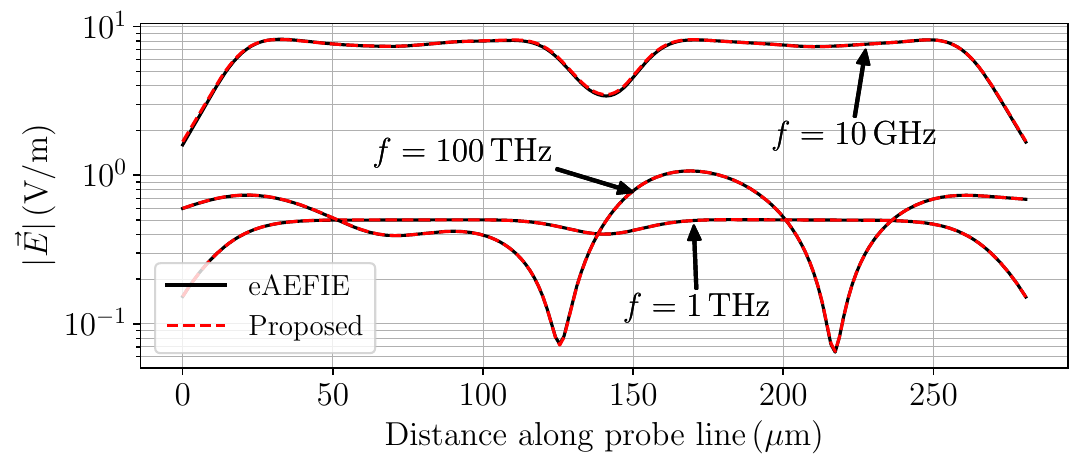}
	\caption{Near-field~$\lvert\Er\rvert$ for the SRR array in \secref{sec:results:srr}.}\label{fig:srr}
\end{figure}

Finally, we consider a~${2 \times 2}$ array of split ring resonators (SRRs), which are of relevance in the design of metamaterials and metasurfaces~\cite{SRRmetaCond}.
Each element has a relative permittivity of~$11$ and an electrical conductivity of~$10^7\,$S/m.
The elements have side length~\SI{2}{\micro\metre}, width~\SI{0.2}{\micro\metre}, and height~\SI{0.1}{\micro\metre}.
The width of the gap is~\SI{0.2}{\micro\metre}.
The structure is meshed with~$3,392$ triangles, and excited with an incident plane wave traveling along the~$-z$ direction, with the electric field polarized in the~$y$ direction.
The geometry is shown in \figref{fig:srrgeom}.
We consider the frequencies~$10\,$GHz,~$1\,$THz and~$100\,$THz.
\figref{fig:srrgeom} shows the electric surface current density for the~$100\,$THz case.
\figref{fig:srr} shows the magnitude of the electric field measured along the probe line shown in \figref{fig:srrgeom}.
The probe line is placed along the~$x$ axis in the~$xz$ plane bisecting the array,~\SI{0.4}{\micro\metre} above it.
Excellent agreement is obtained compared to the eAEFIE for all three frequencies, demonstrating the accuracy of the proposed PIE formulation for realistic structures.

\section{Conclusion}\label{sec:conclusion}
A boundary element formulation based on the electric scalar and magnetic vector potential is proposed for the accurate modeling of lossy objects over wide ranges of frequency and conductivity. Unlike existing potential-based scattering formulations, the proposed method accurately captures the skin effect both at low and high frequencies, and can model both good conductors and lossy dielectrics. The accuracy of the proposed formulation is validated through canonical and realistic numerical examples, and excellent agreement with analytical results and an existing field-based method is observed for at least nine decades of frequency and conductivity.

\clearpage


%

%

%
%

\ifCLASSOPTIONcaptionsoff
  \newpage
\fi


\balance
\bibliographystyle{IEEEtran}
\bibliography{./IEEEabrv,./bibliography}

\end{document}